\documentclass[11pt]{article}

\usepackage[left=1.1in, right=1.1in, top=1.1in, bottom=1.1in]{geometry}

\date{\today}
\newif\ifsiamart
\makeatletter%
\@ifclassloaded{siamart171218}%
  {\siamarttrue}%
  {\siamartfalse}%
\makeatother%

\usepackage{authblk}
\usepackage{silence}
\usepackage[utf8]{inputenc}
\usepackage{xcolor}
\usepackage{graphicx}
\usepackage{bm}
\usepackage{bbm}
\usepackage{amsmath}
\usepackage{amssymb}
\usepackage{stmaryrd}
\usepackage{amsfonts}
\usepackage{mathrsfs}
\usepackage{mathtools}
\usepackage{xparse}
\usepackage{epstopdf}
\usepackage{enumitem}

\ifsiamart\else
\usepackage{amsthm}
\usepackage{hyperref}
\hypersetup{%
    linktocpage=true,
    colorlinks=true,
    citecolor=red,
    linkcolor=blue,
    urlcolor=blue,
}
\usepackage[nameinlink,capitalise]{cleveref}
\newcommand{\email}[1]{\href{mailto:#1}{#1}}

\fi

\definecolor{darkred}{rgb}{.7,0,0}
\definecolor{darkgreen}{rgb}{.1,.7,0}

\usepackage{setspace}

\DeclareMathOperator*{\argmin}{arg\,min}

\newcommand{\placeholder}{\mathord{\color{black!33}\bullet}}%
\newcommand{\ind}{\perp\!\!\!\!\perp}

\newcommand{\normal}{\mathcal{N}}

\newcommand{\real}{\mathbf R}

\renewcommand{\d}{\mathrm d}

\renewcommand{\leq}{\leqslant}

\renewcommand{\le}{\leqslant}

\newcommand{\hmu}{\widehat{\mu}}
\newcommand{\yd}{{y}^\dag}
\newcommand{\hv}{\widehat{v}}

\ifsiamart
\newtheorem{remark}[theorem]{Remark}
\newtheorem{assumption}{Assumption}

\else
\theoremstyle{plain}

\newtheorem{theorem}{Theorem}

\newtheorem{remark}[theorem]{Remark}
\newtheorem{definition}[theorem]{Definition}

\numberwithin{equation}{section}
\crefname{lemma}{Lemma}{Lemmas}
\crefname{remark}{Remark}{Remarks}
\crefname{assumption}{Assumption}{Assumptions}
\crefname{proposition}{Proposition}{Propositions}
\crefname{section}{Section}{Sections}
\crefname{subsection}{Subsection}{Subsections}
\crefname{equation}{}{}
\Crefname{equation}{Equation}{Equations}
\fi

\usepackage{tikz}
\usepackage{tikz-cd}
\usepackage{pgfplotstable}
\pgfplotsset{compat=1.14}
\usetikzlibrary{patterns}
\usetikzlibrary{calc}
\usetikzlibrary{angles}
\usetikzlibrary{quotes}
\usetikzlibrary{external}

\usepackage[url=true,backend=biber,style=trad-abbrv,url=false,doi=false,isbn=false]{biblatex}
\DeclareFieldFormat{volume}{volume \textbf{#1}}
\DeclareFieldFormat[article]{volume}{\textbf{#1}}

\ifsiamart\else
\let\oldparagraph=\paragraph
\renewcommand\paragraph[1]{\oldparagraph{#1.}}
\fi

\newlist{assumpenum}{enumerate}{5}
\setlist[assumpenum]{font={\bfseries}}
\crefname{assumpenumi}{Assumption}{Assumptions}

\definecolor{mypink1}{rgb}{0.858, 0.188, 0.478}
\definecolor{darkgreen}{rgb}{0, 0, 0}
\newcommand{\nnu}{\textcolor{darkgreen}{\mathfrak{\rho}}}
\newcommand{\mdelta}{\textcolor{darkgreen}{\delta}}
\newcommand{\mmu}{\textcolor{darkgreen}{\mathfrak{\mu}}}
\newcommand{\mnu}{\textcolor{darkgreen}{\mathfrak{\nu}}}
\newcommand{\ppi}{\textcolor{darkgreen}{\mathfrak{\pi}}}
\newcommand{\pp}{\textcolor{darkgreen}{\mathfrak{p}}}
\newcommand{\cchi}{\textcolor{darkgreen}{\mathfrak{r}}}
\newcommand{\mapT}{\mathfrak{T}}

\newcommand{\Q}{\mathsf{Q}}
\newcommand{\B}{\mathsf{B}}
\renewcommand{\L}{\mathsf{L}}
\renewcommand{\P}{\mathsf{P}}
\newcommand{\T}{\mathsf{T}}
\newcommand{\G}{\mathsf{G}}
\renewcommand{\S}{\mathsf{S}}
\renewcommand{\S}{\mathsf{S}}

\newcommand{\R}{\mathbb{R}}

\newcommand{\muk}{\mmu^{\rm EK}}

\newcommand{\nuj}{\nnu^{{\rm EK},J}}
\newcommand{\muj}{\mmu^{{\rm EK},J}}
\newcommand{\mujp}{\mmu^{\rm PF}}
\newcommand{\nuk}{\nnu^{\rm EK}}

\renewcommand{\hmu}{\widehat{\mmu}}
\newcommand{\PP}{\mathcal{P}}

\newcommand{\GG}{\mathcal{G}}
\newcommand{\Yd}{Y^{\dagger}}
\renewcommand{\real}{\mathbb R}

\title{%
   {Statistical Accuracy of Approximate Filtering Methods}
}

\author[1]{J. A. Carrillo$^{a,}$}
\author[2]{F. Hoffmann$^{b,}$}
\author[2]{A. M. Stuart $^{c,}$}
\author[3,4]{U. Vaes $^{d,}$}
\affil[ ]{\footnotesize
    $^a$\email{carrillo@maths.ox.ac.uk},
    $^b$\email{franca.hoffmann@caltech.edu},
    $^d$\email{astuart@caltech.edu},
    $^c$\email{urbain.vaes@inria.fr}
}
\affil[1]{\footnotesize Mathematical Institute, University of Oxford, UK}
\affil[2]{\footnotesize Department of Computing and Mathematical Sciences, Caltech, USA}
\affil[3]{\footnotesize MATHERIALS project-team, Inria Paris, France}
\affil[4]{\footnotesize CERMICS, \'Ecole des Ponts, France}

\begin{document}
\maketitle

\section{Introduction}

Estimating the statistics of the state of a dynamical system,
from partial and noisy observations, is both mathematically
challenging and finds wide application. Furthermore, the applications
are of great societal importance, including problems such as
probabilistic weather forecasting~\cite{kalnay2003atmospheric} and prediction of epidemics~\cite{keeling2005networks}.
Particle filters provide a well-founded approach to
the problem, leading to provably accurate approximations of the
statistics~\cite{doucet2001sequential}. However these methods perform
poorly in high dimensions~\cite{MR2459233,snyder2008obstacles}.
In 1994 the idea of ensemble Kalman filtering was introduced~\cite{evensen1994sequential} leading to a methodology that has been
widely adopted in the geophysical sciences~\cite{van2019particle}
and also finds application to quite general inverse problems~\cite{iglesias2013ensemble}. However, ensemble Kalman filters
have defied rigorous analysis of their statistical accuracy, except
in the linear Gaussian setting~\cite{le2009large,mandel2011convergence}.
In this article we describe recent work which takes first steps to
analyze the statistical accuracy of ensemble Kalman filters beyond the
linear Gaussian setting~\cite{carrillo2022ensemble}.
The subject is inherently technical, as it involves the evolution of
probability measures according to a nonlinear and nonautonomous
dynamical system; and the approximation of this evolution.
It can nonetheless be presented in a fairly
accessible fashion, understandable with basic knowledge of
dynamical systems, numerical analysis and probability.
We undertake such a presentation here.

\section{Filtering Problem}

Consider a dynamical system for state $\{v_n\}_{n \in \mathbb{Z}^+}$ evolving in $\real^d$,
partially and noisily observed through data $\{y_n\}_{n \in \mathbb{N}}$ in $\real^K$,
with state and data determined by the following system, holding for $n \in \mathbb{Z}^+$:
        \begin{subequations}
            \begin{alignat*}{3}
                \text{ \bf State: } \qquad
                v_{n+1} &= \Psi(v_{n}) + \xi_{n}\, ,\\
                \text{ \bf Data: } \qquad
                y_{n+1} &= h(v_{n+1}) + \eta_{n+1}\,.
            \end{alignat*}
Here $\mathbb{N}:=\{1,2,\cdots\}$ and $\mathbb{Z}^+=\{0,1,2,\cdots\}.$
We assume that  the initial state of the system is a Gaussian random
variable $v_0 \sim \normal(m_0, C_0)$. Furthermore we assume that $\xi_{n} \sim \normal(0, \Sigma)$ is the
mean-zero noise affecting the state evolution and
$\eta_{n+1} \sim \normal(0, \Gamma)$ is the mean-zero noise entering the data acquisition process.
We assume that the state evolution and data acquisition noise sequences are i.i.d.\ and that the following
independence assumptions hold:
            \[
                v_0 \ind \{\xi_n\}_{n \in \mathbb{Z}^+} \ind \{\eta_{n+1}\}_{n \in \mathbb{Z}^+}.
            \]
        \end{subequations}

The objective of filtering is to determine the probability distribution on $v_n$ given
all the data acquired up to that point. To this end we define, for a given realization
of the data, denoted by a dagger $^\dagger$,
$$\Yd_n  =\{\yd_\ell\}_{\ell=1}^n,
\qquad v_n|\Yd_n \sim \mmu_n.$$
With this notation we can state the objective of filtering more precisely: it is to determine
the \emph{filtering distribution} or \emph{true filter} $\mmu_n$ and update it sequentially in $n$. As we now show $\mmu_n$ evolves according to
a nonautonomous and nonlinear dynamical system on the space of probability measures.

To determine this dynamical system, and various approximations of it that follow,
it is helpful to define
\begin{itemize}
\item $\PP(\real^r):$ all probability measures on $\real^r.$
\item $\GG(\real^r):$ all Gaussian probability measures on $\real^r.$
\end{itemize}
For simplicity we will use the same symbol for measures and their
densities throughout this article.
Before writing the dynamical system for $\mmu_n$ we first determine the
evolution of the measures $\ppi_n, \cchi_n$ defined by
\begin{align*}
v_n &\sim \ppi_n, \qquad (v_n,y_n) \sim \cchi_n.
\end{align*}
These evolve according to
\begin{align*}
\ppi_{n+1}&=\P\ppi_n,\\
\cchi_{n+1}&=\Q\ppi_{n+1}
\end{align*}
Here $\P\colon \mathcal P(\real^d) \to \mathcal P(\real^d)$ is the linear operator
\[
                \P \ppi(u) = \frac{1}{\sqrt{(2\pi)^d \det \Sigma}}\int \exp \left( - \frac{1}{2} |u - \Psi(v)|_{\Sigma}^2 \right) \ppi(v) \, \d v.
            \]
whilst $\Q\colon \mathcal P(\real^d) \to \mathcal P(\real^d \times \real^{K})$, which is also
a linear operator, is determined by
          \[
              \Q \ppi (u, y) = \frac{1}{\sqrt{(2\pi)^K \det \Gamma}}
\exp \left( - \frac{1}{2} \bigl\lvert y - h(u) \bigr\rvert_{\Gamma}^2 \right)
\ppi(u).
          \]
The evolution for $\ppi_n$, which describes the probability of the state~$v_n$, is
determined by a Markov process on $\real^d$ defined via the linear operator $\P.$
The linear operator $\Q$ lifts $\ppi_n$ to the joint space of state and data
$(v_n,y_n).$ It is worth highlighting that, by moving from the evolution of $(v_n,y_n)
\in \real^d \times \real^K$ (finite dimensions) to the evolution of
$(\ppi_n,\cchi_n) \in \PP(\real^d \times \real^K)$ (infinite dimensions) we have converted
a nonlinear stochastic problem into a linear autonomous one.

However, the dynamical system for $\mmu_n$ is nonlinear and nonautonomous.
To determine the dynamical system for $\mmu_n$ we need to introduce a nonlinear
operator on probability measures. Specifically, let~$\B(\placeholder;\yd)\colon \mathcal P(\real^d \times \real^{K})
\to \mathcal P(\real^d)$
describe conditioning of a joint random variable $(v,y)$ on observation $y=\yd$:
          \[
              \B(\nnu;\yd)(u) = \frac{\nnu(u,\yd)}{\int_{\real^d} \nnu(u,\yd) \, \d u}.
          \]
Armed with this we define $v_{n+1}|\Yd_n \sim \hmu_{n+1}$ and observe that
\begin{alignat*}{2}
    \hmu_{n+1}&=\P \mmu_n, \qquad && v_{n+1}|\Yd_n \sim \hmu_{n+1}\\
    \nnu_{n+1}&=\Q\hmu_{n+1}, \qquad &&(v_{n+1},y_{n+1})|\Yd_n \sim \nnu_{n+1}\\
    \mmu_{n+1}&=\B(\nnu_{n+1};\yd_{n+1}), \qquad &&{\rm conditioning}\,{\rm on}\, y_{n+1}=\yd_{n+1}.
\end{alignat*}
Thus we have the nonlinear and nonautonomous dynamical system
\begin{align*}
\mmu_{n+1} & = \B(\Q\P\mmu_n;\yd_{n+1}), \qquad \mu_0=\normal(m_0, C_0).
\end{align*}
This evolution may be thought of in terms of sequential application of
Bayes theorem: $\P\mmu_n$ is prior prediction; $\L(\placeholder;\yd):=\B(\placeholder;\yd) \circ \Q$
maps prior to posterior according to Bayes theorem. This leads
to the following equivalent formulation of the evolution:
\begin{align*}
\mmu_{n+1} & = \L(\P\mmu_n;\yd_{n+1}), \qquad \mu_0=\normal(m_0, C_0).
\end{align*}
The resulting evolution is an infinite dimensional
problem since $\mmu_n$ is a probability density over~$\real^d.$
Approximating this evolution is thus a significant challenge.

\section{Particle Filter}

During the last quarter of the last century, the particle filter became a widely adopted methodology
for solving the filtering problem in moderate dimensions $d$;
the methodology is overviewed in~\cite{doucet2001sequential}.

To describe the basic form of the method we introduce a random
map on the space of probability measures; this map encapsulates Monte
Carlo sampling. To this end let $\Omega$ denote an abstract probability
space (encapsulating the sampling underlying the Monte Carlo method)
and define
$\S^J\colon \PP(\real^r) \times \Omega \to \PP(\real^r)$ to be the empirical
approximation operator:
$$\S^J \mmu = \frac{1}{J} \sum_{j=1}^J \mdelta_{v_j} \, , \qquad v_j \overset{\rm {i.i.d.}}{\sim} \mmu\,.$$
Then for fixed $\omega \in \Omega$ and large $J$, $\S^J$ is a random approximation of the identity on $\PP(\real^r)$, a statement
made precise in the remarks that follow Theorem \ref{t1}.

Since $\S^J \approx I$, we may introduce the approximation
$\mujp_n \approx \mmu_n$ evolving according
to the map
\begin{align*}
\mujp_{n+1} & = \L(\S^J\P\mujp_n;\yd_{n+1}), \qquad \mujp_0=\mmu_0.
\end{align*}
It then follows that, for $n \in \mathbb{N}$,
$$\mujp_{n}=\sum_{j=1}^J w_{n}^{(j)}\mdelta_{\widehat v_{n}^{(j)}}$$
where the particles $\widehat v_{n}^{(j)}$ and weights $w_{n}^{(j)}$
evolve according to
\begin{alignat*}{2}
    \widehat v_{n+1}^{(j)} &= \Psi \bigl(v_n^{(j)}\bigr) + \xi_n^{(j)},\qquad v_{n}^{(j)} \overset{\rm i.i.d.} {\sim} \mujp_{n},\\
            \ell_{n+1}^{(j)} &= \exp \left( - \frac{1}{2} \bigl\lvert \yd_{n+1} - h\bigl(\widehat v_{n+1}^{(j)}\bigr) \bigr\rvert_{\Gamma}^2 \right),\\
w_{n+1}^{(j)}&= \ell_{n+1}^{(j)}\Big/\Bigl(\sum_{m=1}^J \ell_{n+1}^{(m)}\Bigr).
        \end{alignat*}
Here $\xi_{n}^{(j)} \sim \normal(0, \Sigma)$  are independent sequences
of Gaussians each of which is itself i.i.d.\ with respect to both $n$ and $j.$

Systemization of the analysis of this method
may be found in the work of Del Moral~\cite{del1997nonlinear},
with the significant extension to analysis over
long time-intervals in~\cite{del2001stability} and to continuous time
in~\cite{cricsan1999interacting}.
We describe a prototypical theoretical result,
based on the formulation underpinning Rebschini and Van Handel~\cite{MR3375889};
for precise statement of conditions under which the theorem
holds see Law et al~\cite[Theorem 4.5]{MR3363508}.
The theorem deploys a metric $d(\placeholder,\placeholder)$ on random probability measures. This is defined as follows:
$$d(\mmu,\mnu)^2=\sup_{|f| \le 1} \mathbb{E}\bigl|\mmu[f]-\mnu[f]\bigr|^2,$$
where $\mmu[f]=\int f(u) \, \mmu(\d u)$ and similarly for $\mnu;$ the expectation denoted $\mathbb{E}$ is over the randomness used to define random probability measures. The metric reduces to TV distance if applied to  deterministic probability measures.

\begin{theorem}{\cite[Theorem 4.5]{MR3363508}}
\label{t1}
It holds that
\begin{align*}
\sup_{0 \le n \le N} d(\mmu_n,\mujp_{n}) & \le \frac{C(N)}{\sqrt{J}}.
\end{align*}
\end{theorem}

\begin{remark} The following remarks explain and interpret the theorem:

\begin{itemize}
            \item The theorem states that a distance between the particle filter and the true filter is bounded above by
            an error which decreases at the Monte Carlo rate with respect to the number of particles $J$.
            \item The proof follows a typical numerical analysis structure: Consistency + Stability Implies Convergence.
            \item Consistency: $d(\S^J \mmu,\mmu) \le \frac{1}{\sqrt{J}}.$ This shows that $\S^J$ approximates the identity.
            \item Stability: $\P, \L$ are Lipschitz in $d(\placeholder,\placeholder)$.
            \item The theorem assumes upper and lower bounds on the $n$-dependent likelihood arising from sequential application of Bayes theorem; the constant $C$ depends on these bounds. Furthermore $C$ grows exponentially with $N$.
            \end{itemize}
\end{remark}

In practice, especially in high dimensions $d$, the method often suffers from weight collapse~\cite{MR2459233,snyder2008obstacles}.
This refers to the phenomenon where one of the weights $\{w_{n}^{(j)}\}_{j=1}^J$
is close to $1$ and the others are therefore  necessarily close to $0.$
When this happens the method is of little value
because the effective number of particles approximating $\mmu_n$ is $1.$ This leads us to consider
the ensemble Kalman filter which, by design, has equal weights.

\section{Ensemble Kalman Filter}

The ensemble Kalman filter, introduced by Evensen in
\cite{evensen1994sequential}, is overviewed in
the texts~\cite{evensen2009data,evensen2022data}.
We will state the basic particle form of the
algorithm, from~\cite{evensen1994sequential}. We will then
show how it may be derived from a mean-field perspective
highlighted in~\cite{2022arXiv220911371C}.
Finally we describe analysis of this mean-field ensemble Kalman
filter, and its relationship to the true filter~\cite{carrillo2022ensemble}.

To describe the particle form of the method we first introduce key
notation relating to covariances. Specifically we
write covariance under $\nnu \in \mathcal P(\real^d \times \real^{K})$ as:
    \[
        \mathcal {\rm cov}(\nnu) =
        \begin{pmatrix}
            \mathcal C^{vv}(\nnu) & \mathcal C^{vy}(\nnu) \\
            \mathcal C^{vy}(\nnu)^\top & \mathcal C^{yy}(\nnu)
        \end{pmatrix}.
    \]
We use similar notation (${\rm mean}_{\nnu}$) for the mean under $\nnu.$
The ensemble Kalman filter from~\cite{evensen1994sequential} then has the form, for $n \in \mathbb{Z}^+$,
\begin{alignat*}{2}
            \widehat v_{n+1}^{(j)} &= \Psi \bigl(v_n^{(j)}\bigr) + \xi_n^{(j)},\\ 
            \widehat y_{n+1}^{(j)} &= h(\widehat v_{n+1}^{(j)}) + \eta_{n+1}^{(j)}, \\
            v_{n+1}^{(j)} &= \widehat v_{n+1}^{(j)} + \mathcal C^{vy}\bigl(\nuj_{n+1}\bigr) \mathcal C^{yy}\bigl(\nuj_{n+1}\bigr)^{-1} \bigl(\yd_{n+1} - \widehat y_{n+1}^{(j)} \bigr),\\
\nuj_{n+1}&=
\frac{1}{J}\sum_{j=1}^J \delta_{\bigl(\widehat v_{n+1}^{(j)},\widehat y_{n+1}^{(j)}\bigr)}.
\end{alignat*}
Here $\xi_{n}^{(j)} \sim \normal(0, \Sigma)$ i.i.d.\ with respect to both $n$ and $j$
and $\eta_{n}^{(j)} \sim \normal(0, \Gamma)$ i.i.d.\ with respect to both $n$ and $j.$
Furthermore the set of $\{\xi_{n}^{(j)}\}$ is independent of the set of $\{\eta_{n}^{(j)}\}.$

\begin{remark}
\label{rem:emp}
Note that the components of ${\rm cov}(\nuj_{n+1})$ involve empirical learning of the covariances $\Sigma$ and
$\Gamma$ of the $\{\xi_n^{(j)}\}_{j=1}^J$ and $\{\eta_{n+1}^{(j)}\}_{j=1}^J$ respectively. It is possible, and indeed
often desirable, to directly input the matrices $\Sigma$ and $\Gamma$, only using empirical estimation for the covariance based on
the $\{\Psi(v_n^{(j)})\}_{j=1}^J$ and $\{h(\hv_{n+1}^{(j)})\}_{j=1}^J.$
In particular to use
\begin{subequations}
\label{eq:refer}
    \begin{align}
       \mathcal C^{vy}(\nuj_{n+1}) &\approx C^{vh,J},\\
        \mathcal C^{yy}(\nuj_{n+1})& \approx C^{hh,J}+\Gamma
    \end{align}
\end{subequations}
where $C^{vh,J}$ is the empirical cross-correlation between the $\{\Psi(v_n^{(j)})\}_{j=1}^J$ and
the $\{h(\hv_{n+1}^{(j)})\}_{j=1}^J$, whilst $C^{hh,J}$ is the empirical correlation of the
$\{h(\hv_{n+1}^{(j)})\}_{j=1}^J$. The actual definition of $C^{yy}(\nuj_{n+1})$ instead uses an empirical
approximation  $\Gamma_{n+1}$ of $\Gamma$ formed from $\{\eta_{n+1}^{(j)}\}_{j=1}^J.$
\end{remark}

From the particles defined by this algorithm we may define an empirical measure
$$\muj_{n}=
\frac{1}{J}\sum_{j=1}^J \delta_{v_{n}^{(j)}},$$
noting that for $n=0$ we choose $v_{0}^{(j)} \sim \mmu_0$ i.i.d.\,.
Practitioners like this methodology because it assigns equal weights to the
particles and cannot suffer from the weight collapse arising in the particle
filter. Furthermore, in the setting where $\Psi,h$ are both linear, and
$\mmu_n$ is Gaussian, $\muj_{n}$ converges to $\mmu_n$
as $J \to \infty$, at the Monte Carlo rate~\cite{le2009large,mandel2011convergence}.
In this linear setting, with Gaussian noise and initial condition,
the problem can be solved explicitly
by the Kalman filter~\cite{kalman1960new}; this is a nonautonomous and nonlinear
dynamical system for the mean and covariance of the (in this setting) Gaussian $\mmu_n.$
However forming and propagating $d \times d$ covariances, when $d \gg 1$, is impractical;
ensemble methods, on the other hand, operate by computing a $J \times d$ low-rank approximation
of the covariance and can be used when $d \gg 1.$
Thus particle methods are of some practical value even in this linear Gaussian setting,
avoiding the need to work with large covariances when~$d \gg 1.$

We now outline recent new work aimed at determining the statistical accuracy of
the ensemble Kalman filter, beyond the linear Gaussian setting. The mean-field limit
of the particle system, found by letting~$J \to \infty$, is determined by the map
 \begin{alignat*}{2}
            &\widehat v_{n+1} = \Psi (v_n) + \xi_n,\\
            &\widehat y_{n+1} = h(\widehat v_{n+1}) + \eta_{n+1}, \\
            & v_{n+1} = \widehat v_{n+1} + \mathcal C^{vy}\bigl(\nuk_{n+1}\bigr) \mathcal C^{yy}\bigl(\nuk_{n+1}\bigr)^{-1} \bigl(\yd_{n+1} - \widehat y_{n+1} \bigr),\\
            & (\widehat v_{n+1}, \widehat y_{n+1}) \sim \nuk_{n+1}.
        \end{alignat*}
Once again $\xi_{n} \sim \normal(0, \Sigma)$ i.i.d.\ with respect to $n$
and $\eta_{n} \sim \normal(0, \Gamma)$ i.i.d.\ with respect to $n$;
and the set of $\{\xi_{n}\}$ is independent of the set of $\{\eta_{n}\}.$

\begin{remark}
\label{rem:cov}
The covariances $\mathcal C^{vy}$ and $\mathcal C^{yy}$ can be calculated exactly from the covariance of
$\widehat v_{n+1}$ with $h(\widehat v_{n+1})$, from the covariance of $h(\widehat v_{n+1})$
with itself and from $\Gamma$. The formulae are the $J \to \infty$ limit of those given in \eqref{eq:refer}.
\end{remark}

Note that the map $v_n \mapsto v_{n+1}$ is nonlinear, stochastic and nonautonomous in the sense
that it depends on the observed data. But, furthermore, the map depends on the law of $v_n$
since knowledge of this law is required to define $\nuk_{n+1}.$ Thus we now focus on
finding an evolution equation for the law of $v_n$, which we denote  by $\muk_n.$

With this goal we define $\mapT$ as follows:
\begin{alignat*}{2}
\mapT(\placeholder, \placeholder; \nnu,\yd) &\colon
        \real^{d} \times \real^{K} \to \real^d; \\
        (v, y) &\mapsto v + \mathcal C^{vy}(\nnu) \mathcal C^{yy}(\nnu)^{-1} \bigl(\yd - y \bigr),
\end{alignat*}
noting that this map is linear for any given pair $(\nnu,\yd)\in \PP(\real^d\times\real^K)\times\real^K$.
Recalling the notation $F_\sharp {\nu}$ for pushforward of a probability measure $\nu$ under a map~$F$,
we define
\begin{alignat*}{2}
\T(\nnu;\yd) &= \bigl(\mapT(\placeholder, \placeholder; \nnu,\yd)\bigr)_\sharp {\nnu}.
\end{alignat*}
Notice that $\T$, like $\B$, is nonlinear as a map on probability measures.
With the definition of $\T$ it then follows~\cite{2022arXiv220911371C} that
$$\muk_{n+1}  = \T(\Q\P\muk_n;\yd_{n+1}),  \quad \muk_0=\mmu_0.$$
The mean-field map $v_n \mapsto v_{n+1}$ is now well-defined by coupling it to this nonlinear
evolution equation for $\muk_n$, the law of $v_n.$ We refer to it as a \emph{mean-field map}
precisely because of this dependence on its own law.

The key question we wish to address is the relationship between
$\muk_n$ and $\mmu_n$. To focus on this question we write their
evolution equations in parallel:
\begin{align*}
\muk_{n+1} & = \T(\Q\P\muk_n;\yd_{n+1}),\\
\mmu_{n+1} & = \B(\Q\P\mmu_n;\yd_{n+1}).
\end{align*}
From this it is clear that the key is to understand when $\T \approx \B.$
In fact $\T \equiv \B$ on the set of Gaussian measures $\GG(\R^d \times \R^K)$
\cite{2022arXiv220911371C}.
This is related to the fact that, in the Gaussian case, the map $\T$ effects an
exact transport from the joint distribution
of state and data to the conditional distribution of state given data.
To obtain an error bound between $\muk_n$ and $\mmu_n$ it is helpful to define
\begin{align*}
\G &\colon \mathcal{P} \to \mathcal{G},\\
\G\ppi & =\argmin_{\pp \in \mathcal{G}}\,\, d_{\rm KL}(\ppi\|\pp).
\end{align*}
It then follows that
$\G\ppi = \normal({\rm mean}_{\ppi},{\rm cov}_{\ppi})$
\cite[Theorem 4.7]{2018arXiv181006191S}.

\begin{remark}
    \label{eq:trans}
The map $\T$ is not an \emph{optimal} transport, in the sense of transporting one measure into another at minimal cost~\cite{MR2459454}, but it is a transport that is well-adapted to numerical implementation. Links between
transport and data assimilation were pioneered by Sebastian Reich; see~\cite{reich2015probabilistic} and the
references therein.
\end{remark}

When analyzing the particle filter we used  a metric on random
probability measures which reduces to the TV distance in the
non-random case. Here we use a different metric; we have no need to
consider random probability measures, but standard TV alone does not
allow us to control first and second moments. Control of these
moments is natural because of the Gaussian approximations
underlying the use of $\T$ rather than $\B.$ We thus use a
weighted TV metric, with weight $g(v) = 1 + |v|^2$. Specifically we define
    \[
        d_g(\mmu, \mnu) = \sup_{|f| \leq g} \bigl\lvert \mmu[f] - \mnu[f] \bigr\rvert,
    \]
with $\mmu[f]$ and $\mnu[f]$ defined as before. In order to state our theorem about closeness of $\mmu_n$ and $\muk_n$
we introduce the following measure of how close the true filter $\{\mmu_n\}$ is to being Gaussian,
in the lifted space of state and data:

\begin{definition}
We define the closeness between the filtering distribution $\mu_n$, lifted to the joint space
of state and data, and its projection onto Gaussians, over $N$ steps:
    \begin{align*}
\varepsilon:={\rm sup}\,_{0 \le n \le N}\,\, d_g(\G\Q\P\mmu_n,\Q\P\mmu_n).
    \end{align*}
\end{definition}

For precise statement and proof of the following theorem, and in particular details of
the conditions under which it holds, see~\cite{carrillo2022ensemble}.

\begin{theorem}
\label{t4}
Let $\muk_0=\mmu_0$ and assume that
$\|\Psi\|_{L^\infty}, \|h\|_{L^{\infty}}$ and $|h|_{C^{0,1}}$ are bounded by $r$.
Then there is $C:=C(N,r) \in (0,\infty)$ such that
\begin{align*}
\sup_{0 \le n \le N} d_g(\mmu_n,\muk_n) \le C \varepsilon.
\end{align*}
\end{theorem}

\begin{remark} The following remarks explain and interpret the theorem:

\begin{itemize}
            \item The theorem states that a distance between the mean-field ensemble Kalman filter and the true filter is bounded above by a quantity which measures how close the true filter is to being Gaussian.
            This is natural because $\T$ and $\B$ are identical on Gaussians.
            \item As for the particle filter, the proof follows a typical numerical analysis structure: Consistency + Stability Implies Convergence.
            \item Consistency: in~\cite{carrillo2022ensemble}
            a class of problems is identified within which there are sequences of problems
            along which $\varepsilon \to 0$. This demonstrates that there are problem classes
            within which the mean-field ensemble Kalman filter accurately approximates the true filter; at $\varepsilon =0$ a Gaussian problem is obtained and so the result concerns a class of problems in a small neighbourhood of Gaussians.
            \item Stability: $\P, \Q$ are Lipschitz in $d_g(\placeholder,\placeholder)$. Whilst $\B, \T$ are not Lipschitz,
            stability bounds can be proved in $d_g(\placeholder,\placeholder)$, given certain information such as moment
            bounds and lower bounds on covariances.
            \item The theorem assumes upper bounds on the vector fields $\Psi,h$ defining the filtering problems,
            and the constant $C$ depends on these bounds. Furthermore $C$ grows exponentially wth $N$.
            \end{itemize}
\end{remark}

\section{Conclusions}

Ensemble Kalman filters are widely used, yet many open problems remain concerning their
properties. In this article we have concentrated on understanding the sense in which
they approximate the true filtering distribution. This is important for understanding
the sense in which probabilistic forecasts made using ensemble Kalman filters
predict accurate statistics. In the context of weather forecasting, an example
is determining the probability of rain at levels that will cause flooding;
in the context of forecasting epidemics, an example is determining the probability of
an infection peak that will overwhelm health services, in the absence of interventions.
The work described in Theorem \ref{t4} is a first step to build theory in this area,
for non-Gaussian problems. Many challenges lie ahead to develop this theory so that
it applies under more complex and realistic conditions.

\vspace{0.1in}
\paragraph{Acknowledgments}
This article is based on an invited lecture delivered by AMS at ICIAM 2023 in Tokyo.
JAC was supported by the Advanced Grant Nonlocal-CPD (Nonlocal PDEs for Complex Particle Dynamics: Phase Transitions, Patterns and Synchronization) of the European Research Council Executive Agency (ERC) under the European Union’s Horizon 2020 research and innovation programme (grant agreement No. 883363).
JAC was also partially supported by the Engineering and Physical Sciences Research Council (EPSRC) under grants EP/T022132/1 and EP/V051121/1.
FH was supported by start-up funds at the California Institute of
Technology and by NSF CAREER Award 2340762. FH was also supported by the
Deutsche Forschungsgemeinschaft (DFG, German Research Foundation) via
project 390685813 - GZ 2047/1 - HCM.
The work of AMS is supported by a Department of Defense Vannevar Bush Faculty Fellowship,
and by the SciAI Center, funded by the Office of Naval Research (ONR), under Grant Number N00014-23-1-2729.
UV is partially supported by the European Research Council (ERC) under the European Union's Horizon 2020 research and innovation programme (grant agreement No 810367),
and by the Agence Nationale de la Recherche under grants ANR-21-CE40-0006 (SINEQ) and ANR-23-CE40-0027 (IPSO).
The authors are grateful to Eviatar Bach, Ricardo Baptista and Edo Calvello  for helpful suggestions relating to the exposition.

\vspace{0.1in}
\printbibliography

\end{document}